\documentclass[12pt]{amsart}

\usepackage{amssymb,amsthm}

\newtheorem{theorem}{Theorem}
\newtheorem{corollary}{Corollary}
\newtheorem{definition}{Definition}
\newtheorem{lemma}{Lemma}
\newtheorem{problem}{Problem}
\newtheorem{proposition}{Proposition}

\renewenvironment{proof}[1][Proof\ ]{\textit{#1}}{$\hfill{\qed}$}

\newcommand{\ba}{\begin{array}}
\newcommand{\ea}{\end{array}}
\newcommand{\nl}{\newline}
\def\qqquad{\qquad\qquad}
\def\ts{\textstyle}

\def\dim{\mbox{dim}}
\def\deg{\mbox{deg}}
\def\hom{\mbox{Hom}}
\def\End{\mbox{End}}
\def\Ad{\mbox{Ad}}

\def\ker{\mbox{Ker}}

\def\im{\mbox{Im}}
\def\p1{\pi_1}

\def\eq{\Leftrightarrow}

\def\C{{\Bbb C}}
\def\P{{\Bbb P}}
\def\R{{\Bbb R}}

\def\Z{{\Bbb Z}}
\def\V{{\Bbb V}}
\def\G{GL(n,\C)}

\begin{document}

\baselineskip 6mm

\title[Schottky uniformization and vector bundles over Riemann surfaces]
{Schottky uniformization and vector bundles over Riemann surfaces}

\author{Carlos Florentino}

\address{
Department of Mathematics, Instituto Superior T\'ecnico, 
Av. Rovisco Pais, 1049-001, Lisbon, Portugal.}
\email{cfloren@math.ist.utl.pt}

\thanks{   
This work has been partially supported by
FLAD project 50/96 and PRAXIS,
FCT project PCEX/P/MAT/44/96.
}

\date{}

\begin{abstract}
We study a natural map from 
representations of a free group of 
rank $g$ in $GL(n,\C)$, 
to holomorphic 
vector bundles of degree 0 over a
compact Riemann surface $X$ of genus $g$, associated with
a Schottky uniformization of $X$.
Maximally unstable flat bundles
are shown to arise in this way.
We give a necessary and sufficient condition for this map
to be a submersion, when
restricted to representations producing
stable bundles.
Using a generalized version of Riemann's  bilinear relations,
this condition is shown to be true on the subspace 
of unitary Schottky representations.
\end{abstract}

%

\maketitle

\section{Introduction}

Let $X$ be a compact Riemann surface of genus $g$ and
$\mathcal G$ be the sheaf of germs of holomorphic functions from $X$
to $\G$. 
The inclusion $\G\hookrightarrow{\mathcal G}$
(where $\G$ is identified with its constant sheaf on $X$),
defines a map:
\begin{equation}
\label{VT}
{\mathcal V}:H^1(X,\G)\to H^1(X,{\mathcal G}),
\end{equation}
that sends a flat $\G$-bundle into the corresponding
holomorphic vector bundle of rank $n$ over $X$.
Two such flat bundles are said to be {\em analytically equivalent} if
they have the same image under ${\mathcal V}$.
By a slight abuse of terminology, we will also say that a 
holomorphic vector bundle $E$ is {\em flat} if it lies in the image of
${\mathcal V}$; this happens if and only if
every indecomposable component of $E$ has degree 0,
by a classical theorem of Weil \cite{W}.


There is a well know bijection
between the space of flat $\G$-bundles over $X$
and the space of representations 
of the fundamental group of $X$, $\pi_1=\pi_1(X)$,
modulo overall conjugation
$\mathbf{G}_n:=\hom(\pi_1,\G)/\G$;
it is given explicitely by:
\begin{equation}
\label{E.}
E_{\cdot}:\mathbf{G}_n\to H^1(X,\G), \qquad \rho\mapsto 
E_\rho:= \tilde{X} \times_\rho \C^n,
\end{equation}
where $\pi_1$ acts diagonally on the trivial rank $n$ 
vector bundle over the universal cover $\tilde{X}$ of $X$.

The subset 
${\mathcal M}^{st}_n\subset H^1(X,{\mathcal G})$
of stable holomorphic bundles of rank $n$ and degree 0 is a non-singular
quasi-projective algebraic variety of dimension $n^2(g-1)+1$,
as shown by Mumford \cite{M}.
Therefore, if we restrict attention
to the subset 
$\mathbf{G}^{st}_n := \{ \rho:
E_\rho \mbox{ is stable} \}\subset \mathbf{G}_n$
which is a (smooth) complex manifold of twice the dimension of 
${\mathcal M}^{st}_n$ (see \cite{Gu}), we get a map:
\begin{equation}
\label{V}
V:\mathbf{G}^{st}_n \to {\mathcal M}^{st}_n, \qqquad \rho\mapsto 
V(\rho):={\mathcal V}(E_\rho)
\end{equation}
which is surjective, by Weil's theorem 
(a stable bundle is indecomposable)
and is holomorphic (see \S3, below) because of the
universal property of the 
(coarse) moduli space ${\mathcal M}^{st}_n$ (\cite{NS1}).

The well known theorem of Narasimhan and Seshadri \cite{NS2}, 
implies that the restriction of $V$ to the subset
$\mathbf{U}^{st}_n=\{\rho\in\hom(\p1,U(n))/U(n):
\rho \mbox{ is irreducible}\}$
is a diffeomorphism.

Now let us fix a canonical basis of $\pi_1=\pi_1(X)$:
elements $a_1,...,a_g,b_1,...,b_g$ that generate
$\pi_1$, subject to the single relation
$\prod_{i=1}^{g} a_ib_ia_i^{-1}b_i^{-1}=1$.
Let $F_g$ be a free group on $g$ generators $B_1,...,B_g$, and
$p:\pi_1\to F_g$ be the homomorphism given by $p(a_i)=1, p(b_i)=B_i, i=1,...,g$.
Then we can form the exact 
sequence of groups:
\begin{equation}
\label{sequencia}
1 \to N \to \pi_1 \stackrel{p}{\to} F_g \to 1,
\end{equation}
where $N$ is the smallest normal subgroup
of $\pi_1$ containing $a_1,...,a_g$.
A Shottky uniformization,
representing $X$ as a quotient
of a domain of the Riemann sphere 
$\P^1$ by a free subgroup of $PSL(2,\C)$,
gives us a prefered
representation $\sigma_X\in\hom(F_g,PSL(2,\C))$
which is unique up to conjugation (see Thm. \ref{uniformizacao}).
By analogy, we will
say that $\rho\in\hom(\pi_1,\G)$ is a
{\em Schottky representation} if it lies in the image
of the inclusion $i:=p^*:\hom(F_g,\G)\hookrightarrow\hom(\pi_1,\G)$, 
induced from sequence (\ref{sequencia}).
Similarly, a holomorphic vector bundle $E$ is called
a {\em Schottky bundle} (of rank $n$) if 
it is in the image of the composition
\[
\mathbf{S}_n:=\hom(F_g,\G)/\G\stackrel{i}{\hookrightarrow}
\mathbf{G}_n
\stackrel{{\mathcal V}\circ E_{\cdot}}{\to} H^1(X,{\mathcal G}).
\]

In particular, 
a Schottky vector bundle is always flat, and therefore
isomorphic to a direct sum of indecomposable vector
bundles of degree 0.
It is easy to see that all line bundles of degree 0 are Schottky
and that $\mathbf{S}^{st}_n:= \{ \rho:
E_\rho \mbox{ is stable} \}\subset \mathbf{S}_n$ is a 
complex manifold of the same dimension as ${\mathcal M}^{st}_n$
(Prop. \ref{variedades}).

In section 2, we give the first non-trivial 
examples of Schottky vector bundles. 
Let us say that an indecomposable vector bundle $E$ of rank 2
and degree 0 is maximally
unstable, if it is given by a nontrivial
extension $0\to L\to E\to L^{-1}\to 0$, where $L$ is a
square root of the canonical line bundle $K$. 
We will show that:

\begin{theorem}
\label{secundario}
Every maximally unstable indecomposable vector bundle
of rank 2 and degree 0 is Schottky.
\end{theorem}

In section 3 we study the spaces of representations
${\mathbf{G}^{st}_n}$ and ${\mathbf{S}^{st}_n}$, 
and consider the restriction 
\begin{equation}
\label{W}
W:= V|_{\mathbf{S}^{st}_n}: \mathbf{S}^{st}_n\to {\mathcal M}^{st}_n, 
\qqquad \rho\mapsto W(\rho):=V(i(\rho))
\end{equation}
of $V$ to $\mathbf{S}^{st}_n$ which
is again holomorphic (Prop. \ref{holomorfo}).
One can pose the following problem: 
\begin{problem}
Is the above map $W$ surjective, or at least onto
a dense open subset of ${\mathcal M}^{st}_n$ ?
\end{problem}

The answer seems to be unknown at present, to the best of our
knowledge, except for the simplest cases:
rank 1 or genus one\footnote{in a slightly modified version, 
to take into account semistable bundles.},
where it is positive (see \S 6, Appendix).
This problem, which can be called the Schottky uniformization
for vector bundles, may be of interest for the theory of
generalized theta functions, because we can describe
them as holomorphic functions
on $\mathbf{S}^{st}_n$, by pulling back 
the determinant line bundle on the moduli space of semistable bundles
using $W$ (see Beauville \cite{B}).
It can also be useful in describing the K\"ahler metric on the 
moduli space of stable bundles (see Takhtajan and Zograf
\cite{T}, \cite{ZT2}) by analogy
with the fact that the Fuchsian and Schottky uniformizations 
of a compact Riemann surface are related through
a potential for the Weil-Peterson K\"ahler metric
on Teichm\"uller space (\cite{ZT1}).

In section 4, we consider the following period map:
\[
P_{Ad_\rho}:H^0(X,End E_\rho\otimes K) \to H^1(\pi_1,Ad_\rho),
\qquad P_{Ad_\rho}(\phi):=[\gamma\mapsto\textstyle\int_\gamma \phi]
\]
where $\rho\in\mathbf{G}_n$, $Ad_\rho$ denotes the
$\pi_1$-module of $n\times n$ matrices $M$ with the action
$\gamma\cdot M=\rho(\gamma)M\rho(\gamma)^{-1}$,
$H^1(\pi_1,Ad_\rho)$ the first cohomology of $\pi_1$
with values in this module, and $K$ is the
canonical line bundle on $X$. In terms of this period map,
we get a necessary and sufficient condition for $W$ to be locally
invertible.

\begin{theorem}
\label{terciario}
Let $X$ have genus $g\geq 2$.
Then the differential of $W$ at a representation $\rho$,
$dW_\rho:T_\rho\mathbf{S}^{st}_n \to T_{E_\rho}{\mathcal M}^{st}_n$ 
is invertible if and only if 
$H^1(F_g, \Ad_\rho)\cap \im (P_{Ad_\rho}) = 0$.
\end{theorem}

In \S5, we consider the case of representations 
that are both Schottky and unitary. This space
is also important, 
since it is one of the Bohr-Sommerfeld orbits of the
real polarization of the pre-quantum system of 
flat unitary connections on $X$ (see Tyurin \cite{Ty}).
We use a generalization of Riemann's bilinear
relations to prove:

\begin{theorem}
\label{principal}
If $\rho\in\mathbf{S}^{st}_n\cap\mathbf{U}^{st}_n$
then $dW_\rho$ is surjective. In particular,
The image $W(\mathbf{S}^{st}_n)$ 
contains a nonempty open 
set $U\subset{\mathcal M}^{st}_n$ 
(in the complex topology) such that
$W(\mathbf{S}^{st}_n\cap\mathbf{U}^{st}_n)\subset U$.
\end{theorem}

In section 6, we study Schottky bundles in the
easy case $g=1$,
and in the appendix, we consider the case
of line bundles, which have been considered
before from another viewpoint (compare \cite{K}, Ch. VI, \S4).


\section{Schottky uniformization and unstable vector bundles.}

We begin by recalling the classical Schottky uniformization
of compact Riemann surfaces. For details we refer to
\cite{C,AS,Ms,Bs1,Bs2}.
Schottky groups are an important class
of Kleinian groups: groups of M\"obius transformations acting
properly discontinuously on some domain of the Riemann sphere $\P^1$.
A {\em marked Schottky group of genus $g$} is a strictly loxodromic
(this includes hyperbolic)
finitely generated free Kleinian group $\Sigma$ of rank $g$
\cite{Ms}, together with a choice of $g$ generators
$T_1,...,T_g\in PSL(2,\C)$ of $\Sigma$. 
Two marked Schottky groups
$(\Sigma,T_1,...,T_g)$ and $(\Sigma',T_1',...,T_g')$ are 
said to be {\em equivalent} if there
exists a M\"obius transformation $M$ such that $T_i'=MT_iM^{-1}$ for
all $i=1,...,g$.
Thus, the set of equivalence classes of marked Schottky groups of genus $g$
is a subset of $\hom(F_g,PSL(2,\C))/PSL(2,\C)$
(where $F_g$ is
a free group on $g$ generators, and $PSL(2,\C)$ acts by overall
conjugation), called {\em Schottky space of genus $g$}. It is an open set
in $\C^{3g-3}$ \cite{C}.
Let us denote by $\Omega_\Sigma\subset\P^1$
the domain of discontinuity of $\Sigma$
(for $g>1$, its complement is a Cantor set \cite{C},\cite{AS}).

A Schottky group $\Sigma$ gives rise to
a compact Riemann surface $X:=\Omega_\Sigma/\Sigma$.
Every marked Schottky group
$(\Sigma,T_1,...,T_g)$ has a 
(non unique) standard fundamental domain (see \cite{C}),
which is a region $D\subset\P^1$
bounded by smooth closed curves $C_1,...,C_g,C_1',...,C_g'$,
each lying on the outside of all the others,
such that $T_i(C_i)=C_i'$.
If we orient each $C_i$ clockwise and each $C_i'$ counterclockwise,
the canonical holomorphic map
$\Omega_\Sigma \to X$, sends the boundary curves of D
onto smooth non-intersecting simple oriented 
closed curves $\alpha_1,...,\alpha_g$ on $X$.

In this way, we see that a marked Schottky group plus the choice of
a standard fundamental domain determines a
Riemann surface with a distinguished set of 
curves $\{\alpha_1,...,\alpha_g\}$.
Conversely, the
classical retrosection theorem of Koebe
(see \cite{Bs1,AS}), 
states that every compact Riemann surface
arises this way:
{\em
For every compact Riemann surface of genus $g$ with a choice of 
$g$ smooth simple non-intersecting, homologically
independent, oriented closed curves $\alpha_1,...,\alpha_g$,
there exists a
marked Schottky group of genus $g$, $(\Sigma,T_1,...,T_g)$
and a fundamental domain for $\Sigma$ with $2g$ boundary curves
$C_1,...,C_g,C_1',...,C_g'$, 
such that $X=\Omega_\Sigma/\Sigma$ and the map $\Omega_\Sigma\to X$
sends both $C_i$ and $C_i'$ to $\alpha_i$, preserving orientations as above.
The marked Schottky group $(\Sigma,T_1,...,T_g)$ 
satisfying these conditions
is uniquely determined by
$(X,\mathbf{\alpha}_i)$ up to equivalence.}

Recalling that a choice of canonical basis for $\pi_1$ induces a
sequence:
\begin{equation}
1 \to N \to \pi_1 \stackrel{p}{\to} F_g \to 1, 
\end{equation}
this classical theorem can be restated in the following way.
To avoid heavier notaion, we will denote representations 
and their equivalence classes under conjugation
by the same symbols; hopefuly, this should cause no confusion.

\begin{theorem}
\label{uniformizacao}
Let $X$ be a Riemann surface of genus $g\geq 1$ with a
canonical basis for $\pi_1$.
Then there is a unique (up to conjugation) representation 
$\sigma_X\in\!\hom(F_g,\!P\!S\!L(2,\!\C)\!)$
such that $X=\Omega_\Sigma/\Sigma$ where
$\Sigma=\im(\sigma_X)$.
\end{theorem}

\begin{definition}
\label{def-Scho}
A representation $\rho$ of $\pi_1$ in a group $G$
will be called a {\em Schottky representation}, 
if it lies in the image of the
inclusion $i:\hom(F_g,G)\hookrightarrow\hom(\pi_1,G)$.
A holomorphic vector bundle $E$ over $X$ is called a
{\em Schottky vector bundle over $X$} if
$E={\mathcal V}\circ E_\rho$, for some Schottky representation 
$\rho\in \hom(\p1,\G)$ (recall the maps (\ref{VT}, \ref{E.})).
\end{definition}

\noindent 

It is clear that $\rho$ is Schottky if and only if
$\rho(a_i)={\bf 1}$, for all $i=1,...,g$.
Another easy consequence of the definition is that
the tensor product and direct sum of Schottky 
vector bundles is also Schottky.
Later, we will see that every line bundle of degree 0
is Schottky (see \S \ref{linebundles}), 
so the property of being Schottky is preserved under tensoring by
a line bundle of degree 0.

We now give some interesting 
examples\footnote{I thank I. Biswas for 
suggesting to consider these examples.}
of Schottky vector bundles of rank 2, which
are not semistable for $g\geq 2$.
Let $E$ be an extension $0\to L\to E\to L^{-1}\to 0$, where $L$ is a
line bundle of degree $g-1$. For $E$ to be indecomposable, the 
space that classifies these extensions $H^1(X,L^2)=H^0(X,KL^{-2})^*$, 
has to be non zero, 
so $L$ is one of the $2^{2g}$ square roots of the canonical bundle $K$. 
In this case $E$ is unique, and called 
a {\em maximally unstable indecomposable vector bundle}
of rank 2 and degree 0, and denoted here by $E_L$.
We now prove that these bundles are Schottky.

\begin{proof}
{\em of theorem \ref{secundario}}:
Let $L$ be a square root of $K$ (i.e, $L^2=K$). 
First, note that
$E_{L\otimes L_0}=E_L\otimes L_0$, for
any line bundle $L_0$ such that $L_0^2=\mathcal O$; moreover if $E_L$ is
Schottky, so is $E_L\otimes L_0$, by the remarks after definition \ref{def-Scho}. 
Therefore, if the result is true 
for some square root of $K$, then it is also true for any other
square root. If $g=1$, and $\Bbb I$ is the
trivial line bundle, $E_{\Bbb I}$ is Atiyah's bundle ${\Bbb F}_2$ (see section 
\ref{eliptica}), so the result follows from lemma \ref{FibradoAtiyah} below. 
So, assuming $g\geq2$,
consider the following diagram,
induced by the quotient homomorphism $\nu:SL(2,\C)\to PSL(2,\C)$,
whose commutativity is easily established, and 
whose vertical sequences are principal fibrations:
\begin{equation}
\label{diagrama}
\ba{ccccc}
\hom(F_g,\Z_2) & \stackrel{i}{\to} & \hom(\p1,\Z_2) & 
\stackrel{E_\cdot}{\to}
 & H^1(X,\Z_2)\\
\downarrow & & \downarrow & & \downarrow \\
\hom(F_g,SL(2,\C)) & \stackrel{i}{\to} & \hom(\p1,SL(2,\C)) & 
\stackrel{E_\cdot}{\to}
 & H^1(X,SL(2,\C))\\
\downarrow\nu & & \downarrow\nu & & \downarrow\mu \\
\hom(F_g,PSL(2,\C)) & \stackrel{i}{\to} & \hom(\p1,PSL(2,\C))^+ & 
\stackrel{P_\cdot}{\to}
 & H^1(X,PSL(2,\C))^+
\ea
\end{equation}
Here $\hom(\p1,PSL(2,\C))^+$ denotes the connected component 
which contains the image of the inclusion 
$i:\hom(F_g,P\!S\!L(2,\!\C)\!)\hookrightarrow \hom(\p1,P\!S\!L(2,\!\C)\!)$, 
and similarly for $H^1(X,PSL(2,\C))^+$. 
Using the techniques of Gunning (\cite{Gu1}), we see that
the marked Schottky group $\sigma_X\in\hom(F_g,PSL(2,\C))$,
of theorem \ref{uniformizacao}, gives rise to
a flat projective Schottky bundle $P_{i(\sigma_X)}\in 
H^1(X,PSL(2,\C))$, which
by construction is
in the image of the natural map from the space of all projective structures
on $X$ to $H^1(X,PSL(2,\C))$. 
Then, By \cite{Gu1}, th. 2, there is a flat vector
bundle $E_1$ such that $\mu(E_1)=P_{i(\sigma_X)}$ and $E_1$ has divisor order
$g-1$, so that it can be given as an extension
$0\to L_1\to E_1\to L_1^{-1}\to 0$, with $\deg(L_1)=g-1$. 
$E_1$ is indecomposable because of Weil's theorem, thus $L_1^2=K$, and
by uniqueness, $E_1=E_{L_1}$. 
On the other hand, if $\rho\in\hom(F_g,SL(2,\C))$ 
is such that $\nu(\rho)=\sigma_X$,
then $E_{i(\rho)}$ is Schottky by definition, 
and the commutativity of the diagram implies
$\mu(E_{i(\rho)})=P_{i(\sigma_X)}=\mu(E_1)$.
Thus, by exactness, $E_{L_1}$ and $E_{i(\rho)}$, 
differ by tensoring with a $\Z_2$-bundle $L_2\in H^1(X,\Z_2)$
($L_2^2={\mathcal O}$),
so $E_{L_1}$ is Schottky as required.
\end{proof}

Considering vector bundles of the form $E_L\oplus{\Bbb I}^n$,
we get Schottky vector bundles 
which are clearly not semistable, 
for any $g\geq 2$ and any rank $\geq 2$.

\section{Spaces of Schottky representations.}
 
\label{espacos_de_representacoes}

To study the maps $V,W$ of (\ref{V}),(\ref{W}), we start by studying the spaces
of representations $\mathbf{G}_n$ and $\mathbf{S}_n$.
Let $\hom(\pi_1,\G)^o$ and $\hom(F_g,\G)^o$ denote the subsets consisting
of representations $\rho$ having only scalar commutants 
(i.e, if a matrix $M$ commutes with all matrices in the image
of $\rho$, then $M$ is a scalar) and let
\[
\mathbf{G}^{o}_n := \hom(\pi_1,G\!L(n,\!\C)\!)^o/G\!L(n,\!\C)
\qquad
\mathbf{S}^{o}_n := \hom (F_g,G\!L(n,\!\C)\!)^o/G\!L(n,\!\C)
\]
be their quotients under simultaneous conjugation. 
We remark that, in terms of the initial map
${\mathcal V}$ (\ref{VT}), this constitutes no restriction,
since by \cite{Gu},\cite{Gu2} every flat holomorphic
bundle admits a representation with only scalar commutants.
Their {\em tangent spaces} can be given
in terms of cohomology of groups with coefficients in
group modules. 
Denote by $H^k(\Gamma,M)$
the $k$-th cohomology group of $\Gamma$ with coefficients in the 
$\Gamma$-module $M$.
The adjoint representation, together with
a representation $\rho$ of $\pi_1$ (resp. $F_g$) in $\G$, 
endows the space of all $n\times n$ matrices, 
with a $\pi_1$- (resp. $F_g$-) module structure 
($\gamma \cdot M := \rho(\gamma) M \rho(\gamma)^{-1}$
for a matrix $M$), denoted by $\Ad_\rho$.
When $\rho$ has only scalar commutants, it is not difficult to see
that $\dim_\C H^1(F_g,\Ad_\rho) = n^2(g-1)+1$,
which is also the dimension of ${\mathcal M}^{st}_n$.
Note that the inclusion $Z^1(F_g,\Ad_\rho)\hookrightarrow Z^1(\pi_1,\Ad_\rho):
(B_1,...,B_g)\mapsto (0,...,0,B_1,...,B_g)$
induces an inclusion
$H^1(F_g,\Ad_\rho) \hookrightarrow H^1(\pi_1,\Ad_\rho)$, which is
complex linear. 

In \cite{Gu,Gu2} it is shown that $\mathbf{G}^{o}_n$ has the structure
of a complex analytic manifold of dimension $2(n^2(g-1)+1)$ such
that the natural projection $\hom(\pi_1,\G)^o\to\mathbf{G}^{o}_n$
is a complex analytic principal $PGL(n,\C)$-bundle, and
whose tangent space at the equivalence class of $\rho$ is
$H^1(\pi_1,\Ad_\rho)$. $\mathbf{G}^{o}_n$
has also the structure of a complex symplectic
manifold\footnote{$\mathbf{G}^{o}_n$ has also the structure of a
hyperK\"ahler manifold but we will not need this extra structure.},
whose symplectic form $\omega$
is defined by cup product,
an invariant bilinear pairing $B$ on the Lie algebra of $\G$, 
and evaluation on the fundamental homology class $c\in H_2(\pi_1,\C)$,
$\omega : H^1(\pi_1,\Ad_\rho) \times H^1(\pi_1,\Ad_\rho) \stackrel{\cup}{\to}
H^2(\pi_1,\Ad_\rho\otimes\Ad_\rho) \stackrel{B}{\to} H^2(\pi_1,\C) 
\stackrel{\cdot c}{\to} \C$ (see \cite{Go}).
For the case of $\mathbf{S}^{o}_n$ we find:

\begin{proposition}
\label{variedades}
$\mathbf{S}^{o}_n$ is a complex analytic manifold of dimension
$n^2(g-1)+1$ whose tangent space at the equivalence class
of the representation $\rho$ is $H^1(F_g,\Ad_\rho)$. 
Moreover, $\mathbf{S}^{o}_n$ is a
Lagrangian submanifold of $\mathbf{G}^{o}_n$. 
\end{proposition}
\begin{proof}\!\!.
Since $\rho\in\hom(F_g,\G)$ has only scalar commutants, the action
of $PGL(n,\C)$ by conjugation is free and given
$\rho,\rho'\in\hom(F_g,\G)^o$ in the same orbit, there is a unique
$T\in PGL(n,\C)$ such that $\rho'=T\rho T^{-1}$.
The same arguments used in \cite{Gu}, \S9,
prove that $\mathbf{S}^{o}_n=\hom(F_g,\G)^o/PGL(n,\C)$
(the center acts trivially) is a complex manifold,
and that the tangent spaces agree with the required cohomology groups.
Moreover, $\mathbf{S}^{o}_n$ is an analytic
submanifold of $\mathbf{G}^{o}_n$, since the inclusion
$H^1(F_g,\Ad_\rho) \hookrightarrow H^1(\pi_1,\Ad_\rho)$ is complex linear
when $\rho$ is Schottky.
Since $F_g$ is a free group, $H^2(F_g, M)=0$, for any 
$F_g$-module $M$. Therefore, the symplectic form vanishes
on any two tangent vectors to $\mathbf{S}^{o}_n$.
\end{proof}

Observe that $\mathbf{U}^{o}_n=\hom(\p1,U(n))^o/U(n)$
is also a subset of $\mathbf{G}^{o}_n$, 
diffeomorphic to ${\mathcal M}^{st}_n$ by the theorem
of Narasimhan-Seshadri (since a unitary representation
with only scalar commutants is irreducible).
When $\rho$ is a unitary representation of $\pi_1$,
the Lie algebra of $U(n)$ is
again a $\pi_1$-module (denoted $\Ad_\rho^\R$)
via the adjoint representation composed with $\rho$ .
The tangent space of
$\mathbf{U}^{o}_n$ at the equivalence class of the representation
$\rho$ can be identified with $H^1(\pi_1,\Ad_\rho^\R)$ \cite{NS1},
but the inclusion $H^1(\p1,\Ad_\rho^\R)\hookrightarrow H^1(\p1,\Ad_\rho)$
is {\em not} complex linear. Thus, $\mathbf{U}^{o}_n$ sits inside
$\mathbf{G}^{o}_n$ as a {\em real analytic} submanifold
but {\em not} complex analytic,
in contrast to $\mathbf{S}^{o}_n$.

\begin{lemma}
\label{existencia}
For every $n \geq 1$, $g \geq 2$, there are unitary, Schottky
irreducible representations of $\pi_1(X)$ in $\G$. Therefore,
there are stable Schottky vector bundles of any given rank.
\end{lemma}
\begin{proof}\!\!.
We may suppose that $n \geq 2$ (see the appendix). 
Let $\lambda_1, ... , \lambda_n$ be
$n$ distinct complex numbers of modulus 1. Let 
$B_1=diag(\lambda_1,...,\lambda_n)$, and
$B_2$ be the permutation matrix $e_1 \mapsto e_2, ..., e_n \mapsto e_1$, for
a canonical basis $e_1,...,e_n \in {\Bbb C}^n$. It is easy to see that
$B_1, B_2$ form an irreducible set of unitary
matrices (i.e, no subspace of $\C^n$ is preserved by both).
Hence, the representation of $\pi_1$ given by $\rho(a_i)={\bf 1},
\quad i=1,...,g; \quad \rho(b_i)=B_i, \quad i=1,2; \quad
\rho(b_i)={\bf 1}, \quad i=3,...,g$, is unitary, Schottky
and irreducible.
Stability is clear from Narasimhan-Seshadri's theorem.
\end{proof}

\label{familia_universal}
Consider now the subsets $\mathbf{S}^{st}_n$ and $\mathbf{G}^{st}_n$,
of $\mathbf{S}^{o}_n$ and $\mathbf{G}^{o}_n$, respectively,
consisting of conjugacy classes of representations $\rho$
such that $E_\rho$ is stable (since a stable bundle is simple,
$\rho$ has only scalar commutants), and let ${\mathcal E}$ 
be the holomorphic family of vector bundles over $X$
parametrized by the complex manifold $\mathbf{G}^{o}_n$,
whose fiber over $X\times\{\rho\}$ is $E_\rho$.
It is given by the following
$\pi_1$-action on $\mathbf{G}^{o}_n\times\tilde{X}\times\C^n$:
\[
\gamma \cdot (\rho,z,v) = (\rho, \gamma \cdot z, \rho (\gamma) v)
\qqquad \forall \gamma \in \pi_1, (z,v) \in \tilde{X} \times \C^n.
\]

\begin{proposition}
\label{estabilidade'}
Let $X$ be a compact Riemann surface of genus $g\geq2$, and let $n\geq1$.
$\mathbf{S}^{st}_n$ and $\mathbf{G}^{st}_n$
are dense open complex submanifolds
of $\mathbf{S}^{o}_n$ and $\mathbf{G}^{o}_n$, respectively.
\end{proposition}
\begin{proof}\!\!.
In any holomorphic parameter space, 
stable bundles form the complement of an analytic 
subset \cite{NS2}, Thm. 2(B). 
Since $\mathbf{G}^{st}_n$ is non-empty
(irreducible unitary representations
give stable bundles, for $g\geq 2$),
the existence of the family $\mathcal E$ implies that
$\mathbf{G}^{st}_n$ is open and dense
in $\mathbf{G}^{o}_n$.
The same applies to $\mathbf{S}^{st}_n$,
since stable Schottky bundles exist, by lemma \ref{existencia}. 
\end{proof}

\begin{proposition}
\label{holomorfo}
$V: \mathbf{G}^{st}_n \to {\mathcal M}^{st}_n$ and
$W: \mathbf{S}^{st}_n \to {\mathcal M}^{st}_n$ are holomorphic maps. 
\end{proposition}

\begin{proof}\!\!.
This follows from \cite{NS1}, \S 2,
where it is shown that the space of isomorphism classes of
simple vector bundles with degree 0 and rank $n$, ${\mathcal M}^{sim}_n$,
verifies the universal property of a coarse moduli space in the holomorphic
category. In other words, for every holomorphic family 
$\mathcal A$ of simple vector bundles over $X$, 
parametrized by a complex manifold $S$, the ``universal map"
$S\to {\mathcal M}^{st}_n$ sending $s\in S$ to the equivalence class of
${\mathcal A}_s$ is holomorphic. Since our maps $V$ and $W$ are actually 
the ``universal maps" for the holomorphic families 
${\mathcal E}|_{\mathbf{G}^{st}_n}$ and ${\mathcal E}|_{\mathbf{S}^{st}_n}$
constructed above (consisting of stable, hence simple, bundles), and
${\mathcal M}^{st}_n$ is an open subset of ${\mathcal M}^{sim}_n$, we have
the coarse moduli space universal property for ${\mathcal M}^{st}_n$ as 
well.  \end{proof}

We end this section by observing that $\mathbf{S}^{o}_n$ can be
viewed as a natural generalization of Schottky space
considered in \cite{C}.
Recall the canonical map $\nu:\hom(F_g,SL(2,\C))\to\hom(F_g,PSL(2,\C))$.

\begin{proposition}
\label{comutantes_escalares}
Schottky space of genus $g\geq 2$ is contained in 
$\nu(\mathbf{S}^{o}_2)$. More concretely,
if the image of a representation $\sigma:F_g\to PSL(2,\C)$ 
is a Schottky group $\Sigma'$ of genus $g$, 
then $\nu^{-1}(\sigma)$ has only scalar commutants.
\end{proposition}

\begin{proof}\!\!.
By definition of Schottky group, every $T\in\Sigma'=Im(\sigma)$
is loxodromic, and any
two generators of $\Sigma'$ have distinct fixed points. 
Since two non trivial M\"obius transformations
commute if and only if they have the same fixed points, only
the identity commutes with all elements in $\Sigma'$. Therefore, every
representation $\rho\in\nu^{-1}(\sigma)$ has only scalar commutants.
\end{proof}

\section{The period map.}

To compute the differential of $W$,
we now consider a certain period map, and prove
theorem \ref{terciario}.
Given a representation $\sigma:\pi_1\to\G$,
a global holomorphic section of 
$E_\sigma\otimes K$
will be called an {\em $E_\sigma$-differential}.
In terms of a local coordinate $z\in\tilde{X}$, 
an $E_\sigma$-differential
can be viewed as a closed,
holomorphic differential 1-form $\phi=\phi(z)dz$, satisfying
$\phi(\gamma z) \gamma ' (z) = \sigma(\gamma) \phi(z)$,
for all $\gamma \in \pi_1$.
Let $\V_\sigma$ denote the (left) $\pi_1$-module
defined on $\C^n$ by the representation $\sigma$.
For a fixed basepoint $z_0 \in \tilde{X}$ and fixed
$\phi$, simple computations show that the map
$\Phi_{z_0}: \pi_1 \to \V_\sigma$ defined by 
$\gamma \mapsto \int_{z_0}^{\gamma z_0} \phi,$
is a cocycle in $Z^1(\pi_1,\V_\sigma)$, 
(since $\Phi_{z_0}(\gamma_1 \gamma_2)=
\Phi_{z_0}(\gamma_1)+ \gamma_1 \cdot \Phi_{z_0}(\gamma_2)$)
whose equivalence class 
$[\Phi_{z_0}]\in H^1(\pi_1,\V_\sigma)$
does not depend on the basepoint $z_0$ (because
$\Phi_{z_1}(\gamma) = \int_{z_1}^{\gamma z_1} \phi
= (\sigma(\gamma) -1) \int_{z_0}^{z_1} \phi + \Phi_{z_0}(\gamma)$).

\begin{definition}
\label{periodo}
Given a representation $\sigma$, the map:
\[
P_\sigma:H^0(X,E_\sigma\otimes K) \to H^1(\pi_1,\V_\sigma)
\]

defined by $P_\sigma(\phi):=[\Phi_{z_0}]$
is called the {\em period map} associated to $\sigma$.
\end{definition}

If $\rho_1,\rho_2 \in \hom(\pi_1,\G)$, one can ask when are
$E_{\rho_1}$ and $E_{\rho_2}$ analytically equivalent.
Let $z_0\in\tilde{X}$, $p:\tilde{X}\to X$ be the universal cover,
and denote by $f_\phi$ the unique solution  of the differential 
equation $f^{-1}df = p^*\phi, \, f(z_0)=I$,
for a given $\phi \in H^0(X,\End E_{\rho_1} \otimes K)$.

\begin{lemma}
\label{equivalencia}
The following are equivalent: (1) 
$E_{\rho_1} \cong E_{\rho_2}$, \\
(2) there is a holomorphic function 
$f:\tilde{X} \rightarrow \G$ such that \\
$f(\gamma z) = \rho_2 (\gamma) f(z) 
\rho_1(\gamma)^{-1}$
for all $\gamma \in \pi_1, z \in \tilde{X}$, 
and \\
(3) There exists  $\omega \in     H^0(X,\End E_{\rho_1} \otimes K)$
and     $C\in \G$     such that 
$\rho_2 (\gamma) = C f_\omega(\gamma z) \rho_1 (\gamma)  
f_\omega(z)^{-1} C^{-1}$,
for all $\gamma \in \pi_1, z \in \tilde{X}$.
\end{lemma}  

\begin{proof}\!\!. 
To prove $(1)\eq (2)$, note that an isomorphism                    
        between $E_{\rho_1}$ and $E_{\rho_2}$ is a holomorphic global section
 of $E_{\rho_1}^*\otimes E_{\rho_2}=E_{{}^t\rho_1^{-1}\otimes\rho_2}$,  
consisting of invertible matrices. So it corresponds to a holomorphic       
        $f:\tilde{X} \to \G$ such     that         
        $f(\gamma \cdot z) = ({}^t\rho_1^{-1}\otimes\rho_2)(\gamma)     f(z)=
         \rho_2 (\gamma) f(z) \rho_1(\gamma)^{-1}$.\\
To prove $(2)\eq (3)$, let $f$ be as in (2) and put $h(z) :=   
  C^{-1}f(z)$  where $C:=f(z_0)^{-1}$. We obtain $h(z_0)=I$,
        $\rho_2 (\gamma) = C h(\gamma z) \rho_1 (\gamma)  
        h(z)^{-1} C^{-1}$ and                         
        $\omega:=h^{-1}dh$ belongs to $H^0(X,\End E_{\rho_1} \otimes K)$,
since $d(C^{-1}\rho_2(\gamma) C)=0$ is equivalent to
$        ((h^{-1}dh)\circ\gamma) \gamma'   
         = \rho_1(\gamma) (h^{-1}dh) \rho_1(\gamma)^{-1}.$
        Conversely,     if we have (3), then   
        clearly $f=Ch$ verifies (2).
        \end{proof}

For a fixed $\rho\in\hom(\pi_1,\G)$ one can define the following
map 
\[
\ba{lccc}
Q_\rho : & H^0(X, \End E_\rho\otimes K) & \to & \mathbf{G}_n \\
        & \omega  &\mapsto  &
        [f_\omega(\gamma z_0) \rho(\gamma)]
\ea
\]
which, because of the previous lemma,
does not depend on $z_0$, or on the choice of
representative $\rho$ in its conjugacy class.
In fact, using this map, lemma \ref{equivalencia} 
is easily seen to be equivalent to

\begin{lemma}
\label{fibras}
Two bundles $E_{\rho_1}$ and $E_{\rho_2}$ are
analytically equivalent if and only if there is $\omega\in 
H^0(\End E_{\rho_1}\otimes K)$ such that 
$Q_{\rho_1}(\omega)=\rho_2$.
\end{lemma}

If we now restrict to representations $\rho$ with only scalar commutants,
then $Q_\rho$ becomes a holomorphic map, due to the analytic
dependence of $f_\phi$ on $\phi$.

\begin{lemma}
\label{derivada}
{\bf (a)} The differential of $Q_\rho$ at the origin, $d(Q_\rho)_0$,
coincides with $P_{Ad_\rho}$.
\label{ker=im}
{\bf (b)} For any $\rho\in \mathbf{G}^{st}_n$,
$\ker\ dV_\rho=\im\ d(Q_\rho)_0$.
\end{lemma}
\begin{proof}\!\!.
(a) First note that both maps $d(Q_\rho)_0$ and $P_{Ad_\rho}$ are
defined from  \\
$H^0(X,\End E_\rho\otimes K)$ to $H^1(\pi_1,\Ad_\rho)$.
Let $\eta\in H^0(\End E_\rho\otimes K)$
and $t \in \C$. For small $t$, we can expand $f_{t\eta}$ as:
$f_{t\eta}(z) = I + t \int_{z_0}^z \eta + O(t^2)$. 
Let $\rho_t$ denote $Q_\rho(t \eta)$, for brevity.
Discarding second order terms,  we find:
\[
\rho_t(\gamma) = 
        f_{t\eta}(\gamma z_0) \rho(\gamma) =
        \rho(\gamma) + t (\int_{z_0}^{\gamma z} \eta)
        \rho(\gamma) + O(t^2).
\]
The derivative of the curve of representations $\rho_t$ is
given by $\dot{\rho_t}\rho_t^{-1}$,
so the differential at $\phi=0$ in the $\eta$ direction is finally
given by:
\[
\left[ d(Q_\rho)_0 (\eta)\right](\gamma)  =
        \left. \dot{\rho_t}\rho_t^{-1} \right|_{t=0}(\gamma)=
        \lim_{t\to0} \frac{ Q_\rho (t \eta)(\gamma)
        - \rho(\gamma)}{t} \rho(\gamma)^{-1} =
        \int_{z_0}^{\gamma z} \eta \ \nl\indent
\]
(b) For any $\eta\in H^0(X,\End E_\rho\otimes K)$ and $t\in\C$, we
have $V(Q_\rho(t\eta))\cong V(Q_\rho(0))$ (by lemma
\ref{equivalencia}). Letting $t\to0$, we
see that the differential of
$V\circ Q_\rho$, at the origin is zero:
$d(V\circ Q_\rho)_0(\eta)= dV_\rho\circ d(Q_\rho)_0(\eta)=0$. Hence
$\im d(Q_\rho)_0 \subset \ker dV_\rho$.
Conversely, if $\phi\in\ker dV_\rho\subset H^1(\pi_1,\Ad_\rho)$,
then $\phi$ is tangent to the fiber of the map $V$ at $\rho$,
which means tangent to the image of $Q_\rho$ at 0, so
that there is an $\eta\in H^0(X,\End E_\rho\otimes K)$, such that
$\phi=d(Q_\rho)_0(\eta)$.
\end{proof}

\begin{proof}
{\em of theorem \ref{terciario}}:
Since $\mathbf{S}^{st}_n$ and ${\mathcal M}^{st}_n$ 
have the same dimension, $n^2(g-1)+1$
(prop. \ref{variedades}),
we just need to show that 
$dW_\rho$ has trivial kernel.
Since $W$ is the composition
$\mathbf{S}^{st}_n \stackrel{i}{\to} \mathbf{G}^{st}_n 
\stackrel{V}{\to} {\mathcal M}^{st}_n$,
the differential at $\rho$ is the composition:
$H^1(F_g, \Ad_\rho) \stackrel{di}{\to} H^1(\pi_1,\Ad_\rho)
\stackrel{dV_\rho}{\to} T_{E_\rho}{\mathcal M}^{st}_n$,
and so, (by lemma \ref{ker=im}):
\[
\ker dW_\rho = H^1(F_g, \Ad_\rho)
\cap \ker dV_\rho = H^1(F_g, \Ad_\rho)
\cap \im (P_{Ad_\rho}).
\]\end{proof}

\section{Unitary Schottky vector bundles.}

In this section we compute, for unitary Schottky bundles, the
period map using a generalized version of
Riemann's bilinear relations, and prove theorem \ref{principal}.
If $\sigma$ is a unitary representation of
$\pi_1$ in a vector space,
there is a hermitian inner product 
$\langle,\rangle$, which is
invariant under the $\pi_1$ action:
$\langle\gamma \cdot v_1, \gamma \cdot v_2\rangle = 
\langle v_1, v_2\rangle$ for all $\gamma\in\pi_1$ and
vectors $v_1,v_2$ and so, we can define a ``global''
hermitian inner product on the space of $E_\sigma$-
differentials, as follows:
\[
(\phi_1,\phi_2) := 
        i \int_D \langle h_1(z), h_2(z)\rangle \ dz\wedge\overline{dz},
\]
where $D\in\tilde{X}$ is a fundamental domain for $X=\tilde{X}/\pi_1$, and
$\phi_i=h_i(z)dz$ for $z\in\tilde{X}$.
(positive definiteness follows from $\frac i2 dz\wedge\overline{dz}=
dx\wedge dy$).
\label{H1-pairing}
From the Hermitian pairing $\langle,\rangle$, we can also form a pairing in 
$H^1(\pi_1,\V_\sigma)$, as follows. First allow our cocycles
$\Phi\in H^1(\pi_1,\V_\sigma)$ to be defined by linearity
on the group ring
$\Z[\pi_1]$. Let $R_k=\prod_{j=1}^k a_jb_ja_j^{-1}b_j^{-1}$
($k=1,...,g)$, and put $R=R_g$. 
Using the notation of the Fox calculus, we set:
$\frac{\partial R}{\partial a_k} = R_{k-1}-R_kb_k, \quad
\frac{\partial R}{\partial b_k} = R_{k-1}a_k-R_k$.
There is a natural involution $\#$ in $\Z[\pi_1]$ given by 
$\#:\sum n_j \gamma_j \mapsto \sum n_j \gamma_j^{-1}$, so that, for
example
$\#\frac{\partial R}{\partial a_k} = R_{k-1}^{-1}-b_k^{-1}R_k^{-1}, \quad
\#\frac{\partial R}{\partial b_k} = a_k^{-1}R_{k-1}^{-1}-R_k^{-1}$.

In this notation, the fundamental 2-cycle $c\in H_2(\pi_1,\Z)$,
corresponding to $[X]$, under the isomorphism $H_2(X,\Z)\cong H_2(\pi_1,\Z)$,
is given by
$c := \sum_{k=1}^g [\frac{\partial R}{\partial a_k}|a_k] + 
[\frac{\partial R}{\partial b_k}|b_k]$,
where the bar notation $[a|b]$ ($a,b\in\pi_1$) denotes the equivalence class
of the 2-cycle $(1,a,ab)$ in group homology
(see \cite{Br}, ch. II).
Finally, for $\Phi,\Psi\in Z^1(\pi_1,\V_\sigma)$ define:
\[
\Phi\Cup\Psi := \ts \sum_{k=1}^g 
( -\langle\Phi(\#\frac{\partial R}{\partial a_k}), \Psi(a_k)\rangle + 
\langle\Phi(\#\frac{\partial R}{\partial b_k}),\Psi(b_k)\rangle )
\]
This pairing is well defined in cohomology 
(it depends only on the cohomology
classes $[\Phi],[\Psi]\in H^1(\pi_1,\V_\sigma)$) 
and, by an easy calculation, it is the composition
of the cup product, followed by
contraction in $\V_\sigma$ using $\langle,\rangle$, and by
evaluation on the fundamental 2-cycle $c\in H_2(\pi_1,\Z)$
(compare \cite{Go}, \cite{Gu}):
\[
\Cup : H^1(\pi_1,\V_\sigma) \times H^1(\pi_1,\V_\sigma) \stackrel{\cup}{\to}
H^2(\pi_1,\V_\sigma\otimes\V_\sigma) \stackrel{\langle,\rangle}{\to} H^2(\pi_1,\C) 
\stackrel{\cdot \ c}{\to} \C.
\]

\begin{proposition}
\label{bilinear}
{\em (Bilinear relations for $E_\sigma$-differentials).}
\nl
Let $(\V_\sigma,\langle,\rangle)$ be a unitary representation. Then, 
for all $\phi,\psi\in H^0(X, E_\sigma\otimes K)$, we have:
\[
(\phi,\psi) = i \ \{ P_\sigma(\phi)\Cup P_\sigma(\psi) \}.
\]
\end{proposition}
\begin{proof}\!\!.
Let $\Phi(\gamma)=\int_{z_0}^{\gamma z_0}\phi$ be a
cocycle representative of $P_\sigma(\phi)$ and similarly for $\psi$.
Let $f(z):=\int_{z_0}^z \phi$, so that $\phi =
df$ and  $f(\gamma z) = \gamma \cdot f(z)+ 
\Phi(\gamma)$, for all $\gamma\in\pi_1$. By
Stokes' theorem, we have:
\[
(\phi,\psi) 
  = i \int_{\partial D} \langle f(z),\psi(z)\rangle\overline{dz}= 
  i (\sum_{k=1}^{4g} \int_{\gamma_k} \langle f(z),\psi(z)\rangle\overline{dz}),
\]
where the curves $\gamma_k$ are the $4g$ sides of the boundary
of the polygon $D\subset\tilde{X}$,
whose vertices can be ordered as
$\{z_0,a_1z_0,a_1b_1z_0,a_1b_1a_1^{-1}z_0\equiv R_1b_1z_0,
R_1z_0,...,R_gz_0\equiv z_0\}$.
Half of the $4g$ sides give (using the notation $f^\gamma = f\circ \gamma$):
\begin{eqnarray*}
& &\int_{R_{k-1}z_0}^{R_{k-1}a_kz_0} \langle f,\psi\rangle \overline{dz}+ 
\int_{R_{k-1}a_kb_kz_0}^{R_kb_kz_0} \langle f,\psi\rangle \overline{dz}= \\
& & = \int_{z_0}^{a_kz_0} \langle f^{R_{k-1}},\psi^{R_{k-1}}\rangle
    \overline{{R_{k-1}}'(z)} \overline{dz}
    - \int_{z_0}^{a_kz_0} 
    \langle f^{R_kb_k},\psi^{R_kb_k}\rangle \ 
    \overline{(R_kb_k)'(z)} \overline{dz} = \\
& & = \int_{z_0}^{a_kz_0}  \left[ 
    \langle R_{k-1}\cdot f+\Phi(R_{k-1}),R_{k-1}\cdot \psi\rangle -
    \langle R_kb_k \cdot f + \Phi(R_kb_k), 
    R_kb_k \cdot \psi\rangle \right] \overline{dz}=\\
& & = \int_{z_0}^{a_kz_0}  \left[ 
    \langle\Phi(R_{k-1}),R_{k-1}\cdot \psi\rangle -
    \langle\Phi(R_kb_k), R_kb_k \cdot \psi\rangle 
    \right] \overline{dz}=\\
& & = \langle\Phi(R_{k-1}),R_{k-1}\cdot \Psi(a_k)\rangle -
    \langle\Phi(R_kb_k), R_kb_k \cdot \Psi(a_k)\rangle=\\
& & = - \langle\Phi(R_{k-1}^{-1}),\Psi(a_k)\rangle+\langle\Phi(b_k^{-1}R_k^{-1}),\Psi(a_k)\rangle
    = - \langle\Phi(\#\frac{\partial R}{\partial a_k}), \Psi(a_k)\rangle.
\end{eqnarray*}
A similar computation for the remaining $2g$ sides gives the desired formula.
\end{proof}

We observe that, in the special case $\V_\sigma=\C$ with 
the trivial action of $\pi_1$, and the usual
inner product $\langle z_1,z_2\rangle=z_1\overline{z_2}$, $z_1,z_2\in\C$,
this proposition reduces to the
classical Riemann bilinear relations, because it
says that for any
$\phi,\psi\in H^0(X,K)$, we have (here all $R_k$ are trivial):
\[
\ts\int_X \phi\overline{\psi} \; dz\wedge\overline{dz} 
=-i(\phi,\psi)=\Phi\Cup\Psi=
\sum_{k=1}^g 
(\Phi(a_k) \overline{\Psi(b_k)} - 
\Phi(b_k) \overline{\Psi(a_k)}). 
\]

Let us now consider the important case of a representation which is both
Schottky and unitary. Recall that
$P_{Ad_\rho}:H^0(X,End E_\rho\otimes K) \to H^1(\pi_1,\Ad_\rho)$.


\begin{proposition}
\label{corolario}
If $\rho\in\mathbf{S}_n \cap \mathbf{U}_n$ 
then $H^1(F_g, Ad_\rho)\cap \im (P_{Ad_\rho})=0$.
\end{proposition}
\begin{proof}\!\!.
Since $F_g$ is a free group, $H^2(F_g, M)=0$, for any 
$F_g$-module $M$. Therefore, if $P_{Ad_\rho}(\phi) \in H^1(F_g, Ad_\rho)$,
then $P_{Ad_\rho}(\phi) \Cup P_{Ad_\rho}(\phi)=0$
and by proposition \ref{bilinear},
$(\phi,\phi) = 0$. Therefore, $\phi=0$.
\end{proof}

Appying this proposition to representations producing stable
bundles, we can now prove theorem \ref{principal}:

\begin{proof}
{\em of theorem \ref{principal}}:
When $\rho \in \mathbf{S}^{st}_n\cap\mathbf{U}^{st}_n$,
prop. \ref{corolario} and Thm. \ref{terciario}, togheter imply that
$dW_\rho$ is invertible. 
The set $R := \{ \rho \in \mathbf{S}^{st}_n : det (dW_\rho) = 0\}$
is a closed analytic subset of $\mathbf{S}^{st}_n$ which is
not the whole set. Since $W$ is
a holomorphic map between the complex manifolds
$\mathbf{S}^{st}_n$ and ${\mathcal M}^{st}_n$,
on the complement $\mathbf{S}^{st}_n\smallsetminus R$, 
$W$ is a local diffeomorphism,
and therefore, it is an open map. So, we can take
$U = W(\mathbf{S}^{st}_n\smallsetminus R)$. 
\end{proof}

\section{Schottky bundles over an elliptic curve.}
\label{eliptica}

Let $X$ be a compact Riemann surface of genus 1, and ${\Bbb I}$
the trivial line bundle over $X$.
Atiyah \cite{A1} proved the following:

\begin{theorem}
\label{Atiyah}
{\bf (a)} 
For any $n\geq1$, there is a unique indecomposable vector bundle of rank 
$n$ and degree 0 over $X$ denoted ${\Bbb F}_n$, such that
$\dim H^0(X,{\Bbb F}_n)=1$. Moreover, (for $n>1$) ${\Bbb F}_n$ is the
unique nontrivial extension
$0\to {\Bbb I} \to {\Bbb F}_n \to {\Bbb F}_{n-1} \to 0$.
{\bf (b)} 
Every indecomposable vector bundle $E$, of rank $n$ and degree 0 over $X$,
is isomorphic to ${\Bbb F}_n \otimes \det E$.
\end{theorem}

\begin{lemma}
\label{FibradoAtiyah}
For every $n\geq 1$, the bundle ${\Bbb F}_n$ is Schottky.
\end{lemma}
\begin{proof}\!\!.
Write $X=\C/\langle a,b|\ ab=ba\rangle$, where 
$a,b\in\pi_1(X)$ act by
$a\cdot z=z+1, \ b\cdot z=z+\tau$
and $\im\tau>0$.
Consider the Schottky representation $\rho_n\in Hom(F_g,\G)$ given by
assigning to $b$ the $n\times n$ matrix $N$ whose entries are all zero
except for ones on the principal
diagonal and on the diagonal above it
($N_{i,i}=N_{i,i+1}=1$ are the nonzero entries).
We claim that $E_{\rho_n}={\Bbb F}_n$.
Clearly ${\Bbb F}_1=E_{\rho_1}={\Bbb I}$, so assume 
that $E_{\rho_{n-1}}={\Bbb F}_{n-1}$.
By construction, $E_{\rho_n}$ is
an extension of the form
$0\to {\Bbb I} \to E_{\rho_n} \to E_{\rho_{n-1}}={\Bbb F}_{n-1}\to 0$.
and so, by Thm. \ref{Atiyah}(a),
the lemma can be proved by showing that $\dim H^0(X,E_{\rho_n})=1$.
Sections of $E_{\rho_n}$ over $X$
correspond to holomorphic
functions $s=(s_1,...,s_n):\C\to\C^n$ satisfying
$s(\gamma z)=\rho_n(\gamma)s(z)\quad\forall\gamma\in\pi_1$.
This means that $s_n$ is a constant (being an entire doubly
periodic function), and that the other components of $s$
have 1 as period, and verify
$s_i(z+\tau)=s_i(z)+s_{i+1}(z) \quad \forall i=1,...,n-1$. 
Since an abelian differential
with zero ``$a$-periods'' (in this case $ds_{n-1}$) has to be zero,
we get $s_n=0$, and $s_{n-1}$ constant.
Repeating this argument we get $s_i=0$, for all $i=2,...,n$,
and $s_1$ is a constant and so,  $\dim H^0(X,E_{\rho_n})=1$.
\end{proof}

\begin{theorem}
\label{genus_um}
Every flat vector bundle over a Riemann surface of genus one
is Schottky.
\end{theorem}
\begin{proof}\!\!.
By Weil's theorem,
we may assume that $E$ is indecomposable of degree 0.
Then by Thm. \ref{Atiyah}(b), $E={\Bbb F}_n \otimes \det E$.
Since ${\Bbb F}_n$ is Schottky and $\det E$ is a line bundle, we conclude
that $E$ is also Schottky (see the remarks after definition \ref{def-Scho}). 
\end{proof}

Over a Riemann surface of genus one,
there are no {\em stable} vector bundles of degree 0 
and rank $n>1$ (see Tu, \cite{Tu}), so the map $V$ does not make
sense in this case, but we can construct an analog, considering the
moduli space of {\em semistable} vector bundles of degree 0,
which is isomorphic to the $n$-th symmetric product of the Jacobian:
${\mathcal M}^{ss}_n=Sym^n(Jac(X))$ (\cite{Tu}).
Similarly, we can consider semisimple
representations:
\[
\mathbf{G}^{ss}_n=\{\rho\in \hom(\pi_1,\G): \rho(\gamma)
\mbox{ is diagonalizable for all } \gamma\in\pi_1\}/\G,
\]

Since $\pi_1$ is commutative, all matrices $\rho(\gamma)$ can be 
simultaneously diagonalized, and so, $\mathbf{G}^{ss}_n$ is the
$n$-th symmetric product of one dimensional representations,
$\mathbf{G}^{ss}_n=Sym^n(\hom(\pi_1,\C^*))$.
It is then easy to verify that $V_n: \mathbf{G}^{ss}_n\to {\mathcal M}^{ss}_n$ 
(sending $\rho$ to $E_\rho$) is the $n$-th symmetric power of
$V_1: \hom(\pi_1,\C^*)\to {\mathcal M}_1$ (\S\ref{linebundles}).
The space of Schottky representations in
$\mathbf{G}^{ss}_n$, will be
$\mathbf{S}^{ss}_n=Sym^n(\hom(F_g,\C^*))$, and
we can describe the map
$W_1: \mathbf{S}^{ss}_n\to {\mathcal M}^{ss}_n$
using the lattice $\Lambda:=\{2\pi in\omega\}\subset H^0(X,K)$
(where $\omega$ is a normalized
differential), as follows:

\begin{proposition}
For a Riemann surface of genus 1, 
the map $W: \mathbf{S}^{ss}_n\to {\mathcal M}^{ss}_n$ is the $n$-th symmetric
power of the $\Lambda$-bundle $W_1: \hom(F_g,\C^*)\to {\mathcal M}_1$, where
$\Lambda$ acts on $\hom (F_g,\C^*)$ as in Prop. \ref{line_bundles}. 
\end{proposition}

\appendix

\section{Schottky line bundles.}
\label{linebundles}

It is known that every line bundle of degree 0 is Schottky
(see the discussion in \cite{K}, Ch. VI, \S 4).
In our setting, this can be easily obtained as follows.
The moduli space of degree 0 holomorphic line bundles ${\mathcal M}_1$
is the Jacobian
variety of $X$, $Jac(X)$, which is a group under tensor product.
$H^1(X;\C^*)=Hom(\pi_1,\C^*)$ is also an abelian group
(under tensor product of representations)
isomorphic to $(\C^*)^{2g}$
upon the choice of generators of $\pi_1$; 
and the map $V=V_1:\mathbf{G}_1=\hom (\pi_1,\C^*) \to Jac(X)$ 
becomes a (non-algebraic) surjective homomorphism
(a degree 0 line bundle is flat).
The holomorphic map $Q_\rho$ gives now 
an explicit action of $H^0(X,K)$ on $Hom(\pi_1,\C^*)$:
\begin{equation}
\label{accao}
\omega\cdot\rho:=Q_\rho(\omega)(\gamma)= e ^{\int_\gamma \omega} \rho(\gamma),
\end{equation}
which represents the Jacobian as $Hom(\pi_1,\C^*)/H^0(X,K)$.
\label{line_Schottky}
Let $\rho_1\in Hom(\pi_1,\C^*)$ represent a line bundle $L$; 
finding a Schottky representation
$\rho_2$ of the same $L$ amounts, by lemma \ref{fibras}
to finding a holomorphic differential $\omega$ with
$\rho_2(a_j)=e^{\int_{\gamma} \omega}\rho_1(a_j)\\ =1$, for all $j$.
Since this equation is always solved with
$\omega=-\sum_{j=1}^g\log(\rho_1(a_j))\ \omega_j$
(for any choice of branches of log,
and where $\omega_1,...,\omega_g$
is a normalized basis
of $H^0(X,K)$, i.e, $\int_{a_i} \omega_j =\delta_{ij}$ and
$\int_{b_i} \omega_j=\Pi_{ij}$ ($i,j=1,...,g$);
$\Pi_{ij}$ is the period matrix, symmetric
with positive definite imaginary part), we see that
{\em every degree 0 line bundle $L$, admits a Schottky representation.}

Finally, to describe all such representations,
consider the lattice
$\Lambda:=\{2\pi i(n_1\omega_1+...+n_g\omega_g): n_1,...,n_g\in\Z\} $
inside $H^0(X,K)$.
Then,
{\em two Schottky representations $\rho_1$ and $\rho_2$ produce the same
holomorphic line bundle if and only if there exists $\omega\in\Lambda$
such that $Q_{\rho_1}(\omega)=\rho_2$.}
To see this, let $E_{\rho_1}=E_{\rho_2}$
and $\omega\in H^0(X,K)$ be the form such that
$Q_{\rho_1}(\omega)=\rho_2$; then 
$e^{\int_{a_j} \omega} =1$ for all $j$.
Writing $\omega=\sum c_i \omega_i$, this implies $c_j\in 2\pi i\Z$,
which means $\omega\in\Lambda$. The converse is immediate, and 
we conclude that:

\begin{proposition}
\label{line_bundles}
The map $W=W_1:\mathbf{S}_1=\hom (F_g,\C^*) \to {\mathcal M}_1=Jac(X)$ 
is a holomorphic principal
$\Lambda$-bundle, under the action (\ref{accao})
of $\Lambda$ on $\mathbf{S}_1$.
\end{proposition}

Let $\bf 1$ denote the trivial representation.
In contrast to the case of stable bundles,
Schottky bundles do not determine a unique representation:

\begin{corollary}
\label{nao_unicidade}
If $E$ is a Schottky vector bundle, then there are 
infinitely many non-conjugate
Schottky representations that give rise to $E$.
\end{corollary}
\begin{proof}\!\!. If $E=E_\rho$ for some $\rho\in Hom(F_g,\G)$, then
$\rho\otimes (\omega\cdot \bf 1)$
is a non-conjugate Schottky representation, for all 
$\omega\in\Lambda\backslash\{0\}$. Moreover
$E_{\rho\otimes (\omega\cdot{\bf 1})}=E_\rho\otimes E_{\omega\cdot{\bf1}}
=E_\rho\otimes{\Bbb I}=E$, 
by Prop. \ref{line_bundles}.
\end{proof}

{\bf Acknowledgements.}
I am deeply grateful to E. Aldrovandi, E. Bifet, 
P. Zograf, for many interesting discussions, 
and especially to my Ph. D. advisor L. Takhtajan, who introduced
me to this problem. I would like to thank 
also A. Beauville, I. Biswas, M.S. Narasimhan, P. Newstead 
and A. Tyurin, for some useful remarks.




\end{document}